\documentclass[12pt]{article}
\usepackage{amsmath,amsfonts,amsthm,amssymb}
\newcommand\R{{\mathbb R}}
\newcommand\Z{{\mathbb Z}}
\renewcommand\S{\mathbb{S}}
\newcommand\supp{\text{ supp } }
 \newtheorem{theorem}{Theorem}
 
\newtheorem{remark}[theorem]{Remark}
 \newtheorem{lemma}[theorem]{Lemma}

\renewcommand\d{\partial}

\renewcommand\Im{\text{\,Im\,}}
\renewcommand\Re{\text{\,Re\,}}
\begin{document}

\title{Uniqueness of the modified Schr\"odinger map in
$H^{3/4+\varepsilon}(\R^2)$}

\author{Jun Kato \\
Department of Mathematics \\
 Kyoto University \\ jkato@math.kyoto-u.ac.jp 
\and Herbert Koch
 \\ Fachbereich Mathematik   \\ Universit\"at Dortmund \\ koch@math.uni-dortmund.de}

\maketitle

\begin{abstract}
We establish the local well-posedness of the modified Schr\"odinger
map in $H^{3/4+\varepsilon}(\R^2)$.
 \end{abstract}

\section{Introduction}
Let $M$ and $N$ be smooth Riemannian manifolds.
Smooth harmonic maps are smooth maps $f : M \to N$ which locally
minimize the energy
\[ \frac12 \int_M |Df|^2 d\mathrm{vol} \]
where $|Df|$ denotes the Hilbert-Schmidt norm of the derivative of
$f$ and where $d\mathrm{vol}$ is the Riemannian measure on $M$.
They provide an interesting tool for constructing submanifolds
given as the image of $f$. The problem is particularly rich if the
dimension of $M$ is $2$. See \cite{0920.58022, MR2094472}   for the intricate relation between
geometry and analysis in that context.

If $M$ is a Minkowski space instead then the Euler-Lagrange equations
 are the
so called wave map equations, which are an intensively studied prototypical
class of nonlinear wave equations. They occur in several models in
 physics.   See  Tataru
 \cite{MR2043751} and Tao \cite{MR1869874} for a recent results on wave maps.

Their Schr\"odinger version, the Schr\"odinger maps, has been
studied only recently.
In the simplest case one is  given a map $s_0 \in \dot H^1(\R^2; \R^3)$ where
$\dot H^1$ denotes the homogeneous space with
$|s_0|=1$ and one searches a map $ s$ which
satisfies
\[ \d_t s = s \times \Delta s, \]
\[ |s| = 1, \]
\[ s(0,x) = s_0(x). \]
This problem occurs as continuous limit of a cubic lattice of classical spins
evolving in the magnetic field created by their closed neighbors (Sulem, Sulem
and Bardos \cite{MR866199}, Sulem and Sulem \cite{MR1696311}, and  Chang, Shatah and Uhlenbeck
\cite{MR1737504}).

To begin with, we introduce the Schr\"odinger map from
$\R \times \R^{2}$ to the unit sphere $\S^{2}$.
For the map $s(t, \cdot):\,\R^{2} \to \S^{2}$,
we always assume that the energy i.e. the Dirichlet integral of the map
is finite.  In that case the limit
\[  \lim_{R\to \infty} \frac1{\pi R^2}   \int_{B_R(0)}  s(t, x)\, dx \]
exists.
We may and do assume  that this  point is the north pole $N$ after a rotation.
Then, we identify the Riemann surface $(\mathbb{C}, g\,dz\,d\overline{z})$  with
$\mathbb{S}^2\setminus N$ by using stereographic projection, where
$g$ is determined by $g(z, \overline{z})=(1+|z|^{2})^{-2}$ through the relation
\begin{equation*}
  \mathbb{C} \ni z \mapsto
  \Bigl( \frac{2\, \mathrm{Re}\, z}{1+|z|^{2}},\
    \frac{2\, \mathrm{Im}\, z}{1+|z|^{2}},\
    \frac{1-|z|^{2}}{1+|z|^{2}} \Bigr)\in \S^{2},
\end{equation*}
and consider the map $z:\, \R\times\R^{2}\to (\mathbb{C}, g\,dz\,d\overline{z})$.

For each $t\in \R$, the energy of the map $z(t):\, \R^{2}\to
(\mathbb{C}, g\,dz\,d\overline{z})$ is defined by
\begin{equation}\label{energy}
 E(z(t)) = \frac12 \int_{\R^2} \frac{|\nabla z(t)|^2}{(1+|z(t)|^2)^2} dx.
\end{equation}
The Euler-Lagrange equation of the energy functional above is given by
$$
  \sum_{j=1}^{2}
  \left( \frac{\d}{\d x_j} - 2\frac{\overline{z}\, \d_{x_j}
      z}{ 1+|z|^2} \right) \frac{\d z}{\d x_j}=0.
$$
There is a simple geometric meaning of the left hand side: It is the most
natural Laplacian and
\[ \nabla_j = \frac{\d}{\d x_j} - 2\frac{\overline{z}\, \d_{x_j}
      z}{ 1+|z|^2} \]
is the pull back covariant derivative by the map $z$ from $\R\times\R^2$ to
$(\mathbb{C}, g\,dz\,d\overline{z})$.
Here, we notice that $\nabla_0$, the pull back covariant derivative along the direction $t$,  has the similar representation as above denoting
$j=0$ and $x_0=t$.
Then, the Schr\"odinger map equation is given as the evolution equation of
the form
\begin{equation}\label{sm}
  \frac{\d z}{\d t} =  i \sum_{j=1}^{2}
  \left( \frac{\d}{\d x_j} - 2\frac{\overline{z}\, \d_{x_j}
      z}{ 1+|z|^2} \right) \frac{\d z}{\d x_j}.
\end{equation}
Here, we notice that the solution of the equation (\ref{sm}) preserves the energy
(\ref{energy}).

The equation (\ref{sm}) is regarded as a nonlinear Schr\"odinger equation with
a derivative nonlinearity, to which the standard energy method cannot be applied to prove the local well-posedness.
Although there are many studies on the local well-posedness of the initial value
problem for such a class of the derivative nonlinear Schr\"odinger equations
(e.g. Kenig-Ponce-Vega \cite{MR1230709}, Hayashi-Ozawa \cite{MR1255899},
Chihara \cite{MR1731461}), they require many derivatives in $L^2$  of
 the initial data. 
The purpose of this paper is to consider the local well-posedness of the
initial value problem for the Schr\"odinger map (\ref{sm}) for low regularity
initial data, which have a derivative nonlinearity of a specific form.

Following the approach of Nahmod,  Stefanov and Uhlenbeck
\cite{MR1929444} we  apply the ``gauge transformations"
\begin{equation}\label{gauge}
  u_j = e^{i\psi} \frac{\d_{x_j}  z}{1+|z|^2},\quad j=1,2
\end{equation}
to derive the system of the nonlinear Schr\"odinger equations on $u_j$'s by
choosing the ``gauge" $\psi$ appropriately, which
is called the modified Schr\"odinger map equation:
\begin{equation} \label{mse}
\begin{split}
i\d_t u_1 + \Delta u_1 =& -2i A \cdot \nabla u_1 +A_0 u_1
  +|A|^2 u_1+ 4i \Im(u_2\overline u_1) u_2, \\
i\d_t u_2 + \Delta u_2 = & -2i A \cdot \nabla u_2 +A_0 u_2
  +|A|^2 u_2 + 4i \Im(u_1\overline u_2) u_1,
\end{split}
\end{equation}
where $A=(A_1[u,u], A_2[u,u])$ and $A_0=A_0[u,u]$ are defined by
\begin{gather}
 A_j[u,v] =  2\, G_j * \Im (u_1 \overline v_2+v_1 \overline u_2),
 \quad j=1,2, \label{Aj}\\
G_{1} (x)=  \frac{1}{2\pi} \frac{x_{2}}{|x|^2}, \quad
G_{2} (x)= -\frac{1}{2\pi} \frac{x_{1}}{|x|^2},\label{Gj}\\
A_0[u,v]  =  2 \sum_{j,k=1}^2  R_j R_k \Re (u_j \overline
v_k + v_j \overline u_k) +  2 \Re  ( u_1  \overline v_1
+ \overline v_2  u_2), \label{A0}
\end{gather}
and $R_j$ is the Riesz transform defined by the Fourier multiplier
$i\xi_j/|\xi|$.
We summarize the derivation of the modified Schr\"odinger map equation in the appendix.
It is useful to observe that $\nabla \cdot A=0$ and hence $ A \cdot \nabla f=
\nabla \cdot (Af)$ for all $f \in C^1(\R^2)$, which is due to the choice of the
``gauge" $\psi$ (see the appendix).

\begin{remark}
It is not completely obvious whether the Schr\"odinger map problem (\ref{sm})
and the modified Schr\"odinger map problem (\ref{mse}) are equivalent.
This has been shown by Nahmod and Kenig in cases which cover the range of
solutions considered in  this paper.
\end{remark}

For the modified Schr\"odinger maps (\ref{mse}),
Nahmod, Stefanov and Uhlenbeck \cite{MR2038118, MR1929444} have proven
existence of unique solutions for initial data in $H^{1+\varepsilon}(\R^2)$.
The existence part has been extended by the first author and
independently by Nahmod and Kenig \cite{NK} to the construction of solutions
for initial
data in $H^{1/2+\varepsilon}(\R^2)$.
In an even more general context on the original Schr\"odinger map (\ref{sm}),
Ding and Wang \cite{MR1957040, MR1877231} have shown existence of solution
if $u_0 \in H^2$ and uniqueness if $u_0 \in H^4$.
We consider the initial value problem for the modified
Schr\"odinger maps.
Our main result is the following.
\begin{theorem}\label{main}  Suppose that $s > 3/4$ and
$u_0 \in H^s(\R^2)$.
Then there exists a unique local in time solution $u$ to the
modified Schr\"odinger map equation (\ref{mse}) which satisfies
\[ u \in C([0,T]; H^s(\R^2))\cap L^2(0,T;
C^{s-1/2-\varepsilon}(\R^2)) \] for all $T>0$ in the
domain of existence with $\varepsilon >0$ sufficiently small.
More precisely, suppose that
$u_0, v_0 \in H^s(\R^2)$ and that
\[ u,v \in C([0,T]; H^s(\R^2)) \cap
L^2(0,T; C^{s-1/2-\varepsilon}(\R^2)).
\]
Then, with $c$ depending on $\Vert u \Vert_{L^\infty H^{1/2}} $ and
$ \Vert v_0 \Vert_{L^\infty H^{1/2}}$,
\begin{equation}\label{diff}
 \Vert u(t)- v(t) \Vert_{H^{-1/2}}  \le  \exp\bigl\{ c
(1+  \Vert u \Vert_{L^4 B^{1/2}_{4+\varepsilon,2}}^2
+ \Vert v \Vert_{L^4 B^{1/2}_{4+\varepsilon,2}}^2)^2 \bigr\}
\Vert u(s)-v(s) \Vert_{H^{-1/2}}.
\end{equation}
\end{theorem}

The Besov spaces $B^s_{p,q}$ will be  discussed in the next section. There we will see that
\[ L^\infty H^s(\R^2) \cap L^2 C^{\sigma}(\R^2)
  \hookrightarrow L^4 B^{1/2}_{4+\varepsilon,2}(\R^2)
\]
for small $s> 3/4$ and $\sigma > 1/4$.

The existence part is essentially due to the first author in \cite{Kato}.
In \cite{Kato}, the existence of at least one solution was shown for the data
in $H^s$ with $s>1/2$ (see Theorem \ref{thm4.1} below).
The energy estimate of the difference of  two solutions in
$L^2$ implies uniqueness only for the data in $H^1$, which is due to
the loss of derivatives in the nonlinearity.
The main ingredient of this paper is to bound  the difference
of the two solutions in the function space  $H^{-1/2}$ 
to overcome this difficulty.
Indeed, the estimate \eqref{diff} enables us to prove the uniqueness of the solution
for the data in $H^s$ with $s>3/4$, which improves the previous result in $H^1$.

Notation: We denote the standard $L^2$ Sobolev space with $s$ derivatives in
$L^2$ by $H^s$.
The Besov spaces $B^s_{p,q}$ and the H\"older spaces $C^s$ are defined below.
If $X$ is a Banach space we denote by $CX$ resp.~$L^pX$ the spaces of
continuous functions resp.~weakly measurable $p$ integrable functions
with values in $X$,  equipped with the obvious norm.
Constants may change from line to line.  We use $\sim$ with its
standard meaning.

\section{Preliminary results}

We shall use  the Besov space $B^{s}_{p,q}$ with $s \in \R$, $2\le p,q \le
\infty$.
Let $\phi\in C^\infty_0(\R^2)$, $0\le \phi \le 1$, $\phi(\xi)=1$ for $|\xi|\le 1$,
$\phi(\xi)=0$ for $|\xi|\ge 5/4$, $\phi_j(\xi) = \phi(2^{-j} \xi)-
\phi(2^{-j+1} \xi)$ if $ j \ge 1$ and $\phi_0(\xi) =\phi(\xi)$.
Note that $\supp \phi_{j} \subset \{ 2^{j-1}\leq |\xi| \leq 5\cdot2^{j}/4 \}$.
We define $S^j$ through the Fourier multiplier $\phi_j$ and
\[ S_j = \sum_{l=0}^j S^l . \]
We define the spaces $B^s_{p,q}$ through the norm
\[ \Vert f \Vert_{B^s_{p,q}}= \left( \sum_{j=0}^\infty
  2^{sqj} \Vert S^j f
  \Vert_{L^p}^q \right)^{1/q} \]
with the obvious modification if $q=\infty$.
The H\"older space $C^{s}$ with $s>0$ is defined through the norm
$$
  \| f \|_{C^{s}} = \| f \|_{C^{[s]}}+ \sum_{|\alpha|=[s]} \sup_{x\neq y}
    \frac{| \partial^{\alpha}f(x)-\partial^{\alpha} f(y) |}{ |x-y|^{\{s\}}},
$$
where $s=[s]+\{s\}$ with $[s]$ integer and $0\le\{s\}<1$.
If $s\notin \Z_{+}$, it is known that $C^{s}=B^{s}_{\infty,\infty}$ with the equivalent quasi-norm.

Existence of somewhat regular solutions was obtained in
 \cite[Theorem 1.1]{Kato}, which shows the
existence of the solution to the modified  Schr\"{o}dinger map
for the data in $H^{s}(\R^{2})$ with $s>1/2$.

\begin{theorem}\label{thm4.1}
Let $u_{0}\in H^{s}({\R}^{2})$ with $s>1/2$.
Then there exists $T>0$ satisfying
\begin{equation}
   \min\bigl\{1,C/\bigl( (1+ \| u_{0} \|_{L^{2}}^{q})
     \| u_{0} \|_{H^{s}}^{q}\bigr)\bigr\}
    \leq T \leq 1,
\end{equation}
and at least one solution
$u \in L^{\infty}(0,T; H^{s}) \cap C_{w}([0,T]; H^{s})$
to (\ref{mse}) such that
\begin{equation}\label{4.0}
   u \in L^{p}(0,T; B^\sigma_{q,2}(\R^2)),
\end{equation}
for all $2\le p,q \le \infty$ and $\sigma\in \R$ with
\[  s-1/p>\sigma \ge 0, \qquad \text{ and } \qquad 1/p=1/2-1/q.
\]
\end{theorem}

More precisely, given initial data  $u_{0}\in H^{s}(\R^{2})$ with
$s>1/2$, and choosing a small $\varepsilon >0$, a solution $u\in C(L^2)$ is
constructed in Theorem \cite[Theorem 1.1]{Kato} which satisfies
\begin{gather}
  \| (1-\Delta)^{(s-1/2-\varepsilon)/2} u \|_{L^{p}_{T}L^{q}_{x}}
     \leq C \| u_{0} \|_{H^{s}(\R^2)},\label{4.1}\\
  \| u \|_{L^{\infty}_{T}H^{s}}
    \leq C(\|u_{0}\|_{H^{s}}) \| u_{0}\|_{H^{s}(\R^2)},
    \label{4.2}
\end{gather}
where $s-1/2-\varepsilon>2/q> 0$, $1/p=1/2-1/q$, and
$C(\|u_{0}\|_{H^{s}})$ denotes the constant which
depends on $\|u_{0}\|_{H^{s}}$.

The estimate (\ref{4.1}) is the a priori estimate for the  solution
of \cite[Theorem 3.1]{Kato}, and (\ref{4.2}) follows from the energy
estimate \cite[Proposition 2.5]{Kato} combined with (\ref{4.1}).
We choose $q$ large and
apply a Sobolev embedding to see that $q=\infty$ is allowed \eqref{4.1}.
In this process we have to change $\varepsilon$ slightly.
More precisely we even obtain \eqref{4.0} with $q=\infty$.

We  interpolate the norm in \eqref{4.1} with $q=\infty$ and the energy
inequality \eqref{4.2} to obtain the full inequality \eqref{4.0} by using
the following lemma.

\begin{lemma}\label{cal5}
Let $0\leq \theta \leq 1$. We suppose $s_{0}$, $s_{1}\in \R$,
$1\leq q_{0},\ q_{1} \leq \infty$, and
$s=(1-\theta)s_{0}+\theta s_{1}$,
$1/q=(1-\theta)/q_{0}+\theta/q_{1}$.
Then, we have
$$
  \| f \|_{B^{s}_{q,2}} \leq \| f \|_{B^{s_{0}}_{q_{0},2}}^{1-\theta}
    \| f \|_{B^{s_{1}}_{q_{1},2}}^{\theta}.
$$
\end{lemma}

\begin{proof} We estimate
\[ \begin{split}
\sum_{j=0}^\infty (2^{sj} \Vert S^j f \Vert_{L^{q}})^2
\le  &  \sum_{j=0}^\infty (2^{s_0j} \Vert S^j f \Vert_{L^{q_0}})^{2(1-\theta)}
( 2^{s_1 j} \Vert S^j f \Vert_{L^{q_1}})^{2\theta}
\\ \le & \Vert f \Vert_{B^{s_0}_{q_0,2}}^{2(1-\theta)}
 \Vert f \Vert_{B^{s_1}_{q_1,2}}^{2\theta}.
\end{split}
\]
\end{proof}

Now the inequality
\begin{equation} \Vert u(t) \Vert_{B^{s-1/p-\varepsilon}_{q,2}} \le
c\Vert u(t) \Vert_{B^{s-\varepsilon}_{2,2}}^{1-\frac2p}
\Vert u(t) \Vert_{B^{s-\frac12-\varepsilon}
_{\infty,2}}^\frac2p
\end{equation}
holds by using the relation $1/p=1/2-1/q$.
Hence, we obtain
\begin{equation}
\Vert u \Vert^p_{L^pB^{s-1/p-\varepsilon}_{q,2}} \le
\int_{0}^{T}
\Vert u(t) \Vert_{B^{s-\varepsilon}_{2,2}}^{p-2}
\Vert u(t) \Vert_{B^{s-\frac12-\varepsilon}
_{\infty,2}}^2 dt
 \le  \Vert u \Vert_{L^{\infty}H^{s-\varepsilon}}^{p-2}
\Vert u \Vert_{L^2 B^{s-\frac12-\varepsilon}_{\infty,2}}^2.
\end{equation}

Here, we notice that the $H^{s}$-valued strong continuity of the solution in time variable is also obtained  once we prove the uniqueness of the solution in the class of solution in Theorem \ref{thm4.1}.
In fact, it is shown in \cite{Kato} that with a constant depending on 
$\Vert u_0 \Vert_{L^2}$,
\begin{equation}\label{cont} 
\left| \frac{d}{dt} \Vert u_j(t) \Vert_{H^{s}}^2\right|
    \leq C \| (1-\Delta)^{(s-1/2-\varepsilon)/2} u_j(t) \|_{L^{q}}^{2}
    \| u_j(t)\|_{H^{s}}^{2}
\end{equation}
holds for smooth solution $u_j$ converging to the solution $u$ of Theorem 
\ref{thm4.1}, which is constructed regularizing the initial data by the mollifier.
Thus in the  limit 
\[  
\limsup_{s\to 0 } 
 \Vert u(s) \Vert_{H^s} \le \Vert u(0) \Vert_{H^s}.
\]
Since the solution is weakly continuous in time, 
$$  
  \| u(0) \|_{H^{s}} \leq \liminf_{s\to 0} \| u(s) \|_{H^{s}}.
$$
These facts imply the solution is strongly continuous at $t=0$.
Since the solution in this  class is unique by our main result, Theorem
\ref{main},  if $s >3/4$, we are able to apply the argument above for all 
$t \in [0,T]$.
Thus, if $s>3/4$ we obtain 
\begin{equation*}  u \in C([0,T];H^s). \end{equation*}

\section{Calculus in Sobolev spaces}
In this section we establish  several estimates of products in Sobolev spaces and
the Besov spaces, which are used in the proof of Theorem \ref{main}.

\begin{lemma}\label{cal}
Let $n=2$ and $q>4$. Then
the following inequalities are true:
\begin{equation}\label{cal1}  \Vert fg \Vert_{H^{1/2}}
\lesssim \Vert f \Vert_{H^{1/2}}
 \Vert g \Vert_{B^{1/2}_{q,2}},
\end{equation}
\begin{equation} \label{cal2}
\Vert fg  \Vert_{B^{1/2}_{q,2}} \lesssim
\Vert f \Vert_{B^{1/2}_{q,2}}
 \Vert g \Vert_{B^{1/2}_{q,2}}.
\end{equation}
\end{lemma}

\begin{proof}
For the proof of the inequalities \eqref{cal1} and \eqref{cal2} we use the paraproduct
decomposition
\begin{equation}\label{para}
  f \cdot g = \Pi_{1}(f,g) + \Pi_{2}(f,g) + \Pi_{1}(g,f),
\end{equation}
where
\begin{gather*}
  \Pi_{1}(f,g) = \sum_{k=2}^{\infty} S^{k}f \cdot S_{k-2}g,\\
  \Pi_{2}(f,g) = \sum_{k=0}^{\infty} S^{k}f \cdot \widetilde{S}^{k} g,
\end{gather*}
and $\widetilde{S}^{k} = \sum_{l=(k-1)\vee 0}^{k+1} S^{l}$.
Here, we denote $a\vee b=\max(a,b)$ for $a,b\geq 0$.
The important point in the decomposition above is that
\begin{gather*}
  \supp \mathcal{F}[S^{k}f \cdot S_{k-2}g] \subset
    \{ \xi;\, 2^{k-3} \leq |\xi| \leq 2^{k+1} \},\\
  \supp \mathcal{F}[S^{k}f \cdot \widetilde{S}^{k} g] \subset
    \{ \xi;\, |\xi| \leq 2^{k+2} \}
\end{gather*}
hold. 
For the estimate of the $B^{s}_{p,2}$-norm of each term on the right hand
side of (\ref{para}) for $s \geq 0$, $1\leq p \leq \infty$, the support conditions above are used as follows:
\begin{align*}
  \| \Pi_{1}(f,g) \|_{B^{s}_{p,2}}^{2}
  & = \sum_{j=0}^{\infty} 2^{2sj} \| S^{j} \Pi_{1}(f,g) \|_{L^{p}}^{2}\\
  & = \sum_{j=0}^{\infty} 2^{2sj} \Bigl\| S^{j} \Bigl(
     \sum_{k=(j-1)\vee 2}^{j+3}
    S^{k}f \cdot S_{k-2}g \Bigr) \Bigr\|_{L^{p}}^{2}\\
  & \lesssim \sum_{k=2}^{\infty} 2^{2sk}
    \| S^{k}f \cdot S_{k-2}g \|_{L^{p}}^{2},
\end{align*}
and similarly
\begin{align*}
  \| \Pi_{2}(f,g) \|_{B^{s}_{p,2}}
  & = \Bigl\{ \sum_{j=0}^{\infty} 2^{2sj} \Bigl\| S^{j} \Bigl(
    \sum_{k=(j-2)\vee 0}^{\infty}
    S^{k}f \cdot \widetilde{S}^{k} g \Bigr)
    \Bigr\|_{L^{p}}^{2} \Bigr\}^{1/2} \\
  & \lesssim \Bigl\{ \sum_{j=0}^{\infty} \Bigl(
    \sum_{k=(j-2)\vee 0}^{\infty} 2^{sj}
    \| S^{k}f \cdot \widetilde{S}^{k} g \|_{L^{p}} \Bigr)^{2} \Bigr\}^{1/2}\\
  & \lesssim \sum_{k=0}^{\infty} 2^{sk}
    \| S^{k}f \cdot \widetilde{S}^{k} g \|_{L^{p}}.
\end{align*}

Then, it is easy to see the estimate
\begin{equation} \label{two}
  \| f g \|_{B^{1/2}_{q,2}} \lesssim \| f \|_{B^{1/2}_{q,2}} \| g \|_{B^{0}_{\infty,1}}
    + \| f \|_{B^{0}_{\infty,1}} \|g \|_{B^{1/2}_{q,2}},
\end{equation}
holds, because
\begin{gather*}
\Vert \Pi_1(f,g) \Vert_{B^{1/2}_{q,2}} 
  \lesssim \Vert g \Vert_{L^\infty}  \Vert f \Vert_{B^{1/2}_{q,2}}
  \lesssim \Vert g \Vert_{B^{0}_{\infty,1}}  \Vert f \Vert_{B^{1/2}_{q,2}},\\
\Vert \Pi_2(f,g) \Vert_{B^{1/2}_{q,2}}
  \lesssim \Vert g \Vert_{B^{0}_{\infty,1}}  \Vert f \Vert_{B^{1/2}_{q,2}}.
\end{gather*}
Inequality  \eqref{cal2}  is an  immediate consequences of the inequality above
since $B^{1/2}_{q,2}\hookrightarrow B^{0}_{\infty,1}$ for $q>4$.

Now we turn to the proof of the inequality (\ref{cal1}).
Since $H^{s}=B^{s}_{2,2}$ with the equivalent norms, we have
\begin{align*}
  \| \Pi_{1}(f,g) \|_{H^{1/2}}
  & \lesssim \Bigl\{ \sum_{k=2}^{\infty} 2^{k}
    \| S^{k}f\|_{L^{2}}^{2} \|S_{k-2}g \|_{L^{\infty}}^{2} \Bigr\}^{1/2}\\
  & \lesssim \| g \|_{L^{\infty}} \Bigl\{ \sum_{k=2}^{\infty} 2^{k}
    \| S^{k} f \|_{L^{2}}^{2} \Bigr\}^{1/2}\\
  & \lesssim \| g \|_{B^{1/2}_{q,2}} \| f \|_{H^{1/2}},
\end{align*}
since $B^{1/2}_{q,2} \hookrightarrow L^{\infty}$ for $q>4$.
For the estimate of the other terms, we set $1/p+1/q=1/2$ with $q>4$, and
we use the Sobolev embedding
\begin{equation}\label{young}
  \| S^{k} f \|_{L^{p}} \lesssim
    2^{2k/q} \| S^{k} f \|_{L^{2}}
\end{equation}
to obtain
\begin{align*}
  \| \Pi_{2}(f,g) \|_{H^{1/2}}
  & \lesssim \sum_{k=0}^{\infty} 2^{k/2}
    \| S^{k}f \|_{L^{p}} \| \widetilde{S}^{k} g \|_{L^{q}}\\
  & \lesssim \sum_{k=0}^{\infty} 2^{k/2}  \| S^{k}f \|_{L^{2}}
    2^{2k/q} \| \widetilde{S}^{k} g \|_{L^{q}}\\
  & \lesssim \| f \|_{H^{1/2}} \| g \|_{B^{2/q}_{q,2}}.
\end{align*}
Note that $ B^{1/2}_{q,2} \hookrightarrow B^{2/q}_{q,2}$ holds for $q>4$.
Finally, we use the inequality (\ref{young}) again to obtain
\begin{align*}
  \| \Pi_{1}(g,f) \|_{H^{1/2}}
  & \lesssim \Bigl\{ \sum_{k=2}^{\infty} 2^{k}
    \| S_{k-2}f\|_{L^{p}}^{2} \|S^{k}g \|_{L^{q}}^{2} \Bigr\}^{1/2}\\
  & \lesssim \Bigl\{ \sum_{k=2}^{\infty} 2^{k}
    \Bigl( \sum_{l=0}^{k-2} \| S^{l} f \|_{L^{p}} \Bigr)^{2}
    \|S^{k}g \|_{L^{q}}^{2} \Bigr\}^{1/2}\\
  & \lesssim \Bigl\{ \sum_{k=2}^{\infty} 2^{k}
    \Bigl( \sum_{l=0}^{k-2} 2^{2l/q} \| S^{l} f \|_{L^{2}} \Bigr)^{2}
    \|S^{k}g \|_{L^{q}}^{2} \Bigr\}^{1/2}\\
  & \lesssim \sum_{l=0}^{\infty}
    \Bigl\{ \sum_{k=l+2}^{\infty} 2^{k} \|S^{k}g \|_{L^{q}}^{2} \Bigr\}^{1/2}
    2^{2l/q} \| S^{l} f \|_{L^{2}}\\
  & = \| g \|_{B^{1/2}_{q,2}} \sum_{l=0}^{\infty} 2^{-(1/2-2/q)l}
    2^{l/2} \| S^{l} f \|_{L^{2}}\\
  & \lesssim \| g \|_{B^{1/2}_{q,2}} \| f \|_{H^{1/2}}
\end{align*}
for $q>4$.
This completes the proof of the inequality (\ref{cal1}).
\end{proof}

\begin{lemma} The following inequalities are true
\begin{equation}\label{cal3}
 \Vert f g \Vert_{H^{-1/2}} \lesssim \Vert f \Vert_{H^{-1/2}} \Vert g
\Vert_{B^{1/2}_{q,2}}
\end{equation}
and
\begin{equation}\label{cal4}
\Vert f\nabla g \Vert_{H^{-1/2}} \lesssim \Vert f \Vert_{H^{1/2}}
  \Vert g \Vert_{B^{1/2}_{q,2}}.
\end{equation}
\end{lemma}
\begin{proof}
The first inequality follows by duality (consider $f \to fg$ with fixed $g$) from  the previous lemma.
For the second inequality we use $ f \nabla g = \nabla (fg) - (\nabla f) g$
and estimate
\[\begin{split}
 \Vert  f \nabla g \Vert_{H^{-1/2}} \le & \Vert \nabla (fg) \Vert_{H^{-1/2}} +
 \Vert (\nabla f) g \Vert_{H^{-1/2}} \\
\lesssim  &  \Vert fg \Vert_{H^{1/2}} + \Vert \nabla f \Vert_{H^{1/2}}
\Vert g \Vert_{B^{1/2}_{q,2}}
\\ \lesssim & \Vert f \Vert_{H^{1/2}} \Vert g \Vert_{B^{1/2}_{q,2}}.
\end{split}
\]
\end{proof}

In what follows, we collect several estimates on the vector field
$A[f,g]=(A_1[f,g], A_2[f,g])$ defined by \eqref{Aj}, \eqref{Gj} by using the
estimates above.
It is easy to see that each $A_j$ is real, and $\nabla \cdot A =0$ holds.

\begin{lemma}
The vector field  $A$ defined by \eqref{Aj}, \eqref{Gj} satisfies
\begin{gather}
\Vert \nabla A[f,g] \Vert_{\infty}
\lesssim \Vert f \Vert_{B^{1/2}_{q,2}} \Vert g \Vert_{
B^{1/2}_{q,2}}
+  \Vert f \Vert_{L^2} \Vert g \Vert_{L^2}, \label{A}\\
\Vert A[f,g] \Vert_{B^{1/2}_{q,2}}
  \lesssim  \Vert f \Vert_{B^{1/2}_{q,2}} \Vert g \Vert_{B^{1/2}_{q,2}}
  + \Vert f \Vert_{L^2} \Vert g \Vert_{L^2}, \label{Abesov}
\end{gather}
for $q>4$.
\end{lemma}

\begin{proof}
The inequality \eqref{A} is derived as follows.
\begin{equation}
\begin{split}
\Vert \nabla A[f,g] \Vert_{L^{\infty}}
\lesssim & \Vert fg \Vert_{C^\varepsilon} +
\Vert fg \Vert_{L^1}
\\ \lesssim & \Vert f \Vert_{C^\varepsilon} \Vert g
\Vert_{C^{\varepsilon}}
+ \Vert f \Vert_{L^2} \Vert g \Vert_{L^2}
\\ \lesssim  &
 \Vert f \Vert_{B^{1/2}_{q,2}}
\Vert g \Vert_{B^{1/2}_{q,2}} + \Vert f \Vert_{L^2} \Vert g
\Vert_{L^2},
\end{split}
\end{equation}
where we used elliptic regularity and the embedding
$B^{1/2}_{q,2}(\R^2) \subset C^\varepsilon(\R^2)$.

To derive the inequality \eqref{Abesov} we decompose $A[f,g]$ into the high
frequency part and the low frequency part,
\[ A[f,g]= S_0 A[f,g] + (I-S_0) A[f,g], \]
where $S_0$ is the operator defined in Section 2, which is the Fourier multiplier
$\phi$.
Then, the high frequency part is easily estimated as
\begin{align*}
  \| (I-S_0) A[f,g] \|_{B^{1/2}_{q,2}}
  & \lesssim \| fg \|_{B^{-1/2}_{q,2}}\\
  & \lesssim \| fg \|_{B^{1/2}_{q,2}}\\
  & \lesssim \| f \|_{B^{1/2}_{q,2}} \| g \|_{B^{1/2}_{q,2}}.
\end{align*}
Here the first inequality follows from the elliptic regularity (it is not hard to see
that there no difficulty from the low frequency part), and the third from
\eqref{cal1}.

To estimate the low frequency part, we first observe that
\begin{align*}
  S_0 A[f,g] & =  F^{-1}[\phi]\ast A[f,g]\\
  & \sim F^{-1}[\phi]\ast K \ast (fg),
\end{align*}
where $F^{-1}$ denotes the inverse Fourier transform and $K$ is the smooth
homogeneous function of degree $-1$.
Then, it is easy to see $\Phi \equiv F^{-1}[\phi]\ast K \in L^{r}(\R^2)$ for
$2<r<\infty$ by using the Hardy-Littlewood-Sobolev inequality.
Thus, we obtain by Young's inequality
\begin{align*}
  \| S_0 A[f,g] \|_{B^{1/2}_{q,2}} & \lesssim  \Vert S_0 A[f,g]  \Vert_{L^q} \\
& = \| \Phi \ast (fg) \|_{L^q}\\
  & \lesssim \| \Phi \|_{L^{q}} \| f \|_{L^{2}} \| g \|_{L^{2}},
\end{align*}
where we also use the embedding $B^{1/2+\varepsilon}_{q,q} \hookrightarrow
B^{1/2}_{q,2}$, the second inequality follows from the fact that we are concerned
with the low frequency part, and the third from the fact
$B^{0}_{q,q}=L^q$ with the equivalent norm.
This completes the proof of the inequality \eqref{Abesov}.
\end{proof}

\begin{lemma}\label{Abi}
Let $n=2$ and $q>4$. Then the following inequality holds.
\begin{equation}
 \| A [f,g] \nabla h  \|_{H^{-1/2}} \lesssim  \bigl( \| g \|_{H^{1/2}} \| h \|_{H^{1/2}}
    + \| g \|_{B^{1/2}_{q,2}} \| h \|_{B^{1/2}_{q,2}} \bigr)
   \| f \|_{H^{-1/2}},
\end{equation}
where $A[f,g]$ is the vector field defined by \eqref{Aj}, \eqref{Gj}.
\end{lemma}

\begin{proof}
We use   that $A$ is divergence free and hence
\[
\Vert A[f,g] \nabla h \Vert_{H^{-1/2}} = \Vert \nabla (A[f,g] h)
\Vert_{H^{-1/2}} \lesssim  \Vert A[f,g]h \Vert_{H^{1/2}}.
\]
As in the proof of the previous lemma, we decompose $A[f,g]$ into the high
frequency part and the low frequency part,
\[ A[f,g]= S_0 A[f,g] + (I-S_0) A[f,g]. \]
Then, the high frequency part is easily estimated as follows.
\begin{equation}
\begin{split}
\Vert (I-S_0)  A[f,g]\cdot h \Vert_{H^{1/2}}
\lesssim   &  \Vert  h \Vert_{B^{1/2}_{q,2}}
\Vert (I-S_0)A[f,g] \Vert_{H^{1/2}} \\
\lesssim &  \Vert  h \Vert_{B^{1/2}_{q,2}}
\Vert fg \Vert_{H^{-1/2}} \\
\lesssim &
\Vert  h \Vert_{B^{1/2}_{q,2}}
\Vert g \Vert_{B^{1/2}_{q,2}}
\Vert f \Vert_{H^{-1/2}}.
\end{split}
\end{equation}
Here the first inequality follows from \eqref{cal1}, the second from
elliptic regularity, and the third inequality from \eqref{cal3}.

For the estimate on the low frequency part, we observe that a direct
calculation shows that
\[ \Vert uv \Vert_{L^2} \le \Vert u \Vert_{L^\infty} \Vert v \Vert_{L^2} \]
and
\[ \Vert \nabla (uv) \Vert_{L^2} \le
(\Vert u \Vert_{L^\infty} + \Vert \nabla u \Vert_{L^{\infty}}) \Vert v
\Vert_{H^1},
\]
hence by interpolation
\[ \Vert S_0 A[f,g] \cdot h \Vert_{H^{1/2}}
\le 2 \Vert S_0 A[f,g] \Vert_{W^{1,\infty}} \Vert h \Vert_{H^{1/2}}.
\]
To complete the proof we have to bound $S_0A[f,g]$ and its gradient.
By translation invariance it suffices to do this at the origin.
The argument for
$S_0A[f,g]$ and for its gradient is the same.
Now using the notation in the proof of the previous lemma, we observe that
\begin{align*}
  S_0 A[f,g](0) & \sim  \{ \Phi \ast  (fg) \} (0)\\
  & = \int_{\R^2} \Phi(y) f(y) g(y) dx.
\end{align*}
Thus,
\begin{align*}
  | S_0 A[f,g](0) |
  & \lesssim \Bigl| \int_{\R^{2}} \Phi(y) f(y) g(y) dy \Bigr| \\
  & \lesssim \| \Phi\, g \|_{H^{1/2}} \| f \|_{H^{-1/2}}\\
  & \lesssim \| \Phi \|_{B^{1/2}_{q,2}} \| g \|_{H^{1/2}} \| f \|_{H^{-1/2}}.
\end{align*}
Finally, we notice that
$$
  \| \Phi \|_{B^{1/2}_{q,2}} \lesssim \| \Phi \|_{B^{1/2+\varepsilon}_{q,q}}
    \lesssim \| \Phi \|_{B^0_{q,q}} \sim \| \Phi \|_{L^{q}} < \infty,
$$
since $\Phi$ is supported in the low frequency part in the Fourier space, and
$\Phi \in L^{r}(\R^{2})$ for $2<r<\infty$.
This completes the proof.
\end{proof}

\section{The energy estimates}
In this section we derive elementary $L^2$ inequalities for Schr\"odinger
equations  with a drift term. We first consider
\[ i\partial_t u + \Delta u + iv\cdot \nabla  u = F. \]
\begin{lemma}\label{l2e}  Suppose that $v$ is real valued and
\[ \int_0^T \Vert \nabla v(\tau)  \Vert_{L^{\infty}}d\tau < \infty,
\]
where $\| \nabla v \|_{L^{\infty}}=\|  |\nabla v|_{HS}
\|_{L^{\infty}}$  with $ | \cdot |_{HS}$ the Hilbert-Schmidt
norm, more precisely the Euclidean length of the vector in the case here.
Then the following inequality holds for  $0\le s \le 1$
\[
 \Vert u(t) \Vert_{H^s(\R)}
\le e^{ 4 \int_0^t \Vert \nabla v \Vert_{L^{\infty}} d\tau }
\Bigl(\Vert u(0)
\Vert_{H^s(\R)}  + \int_0^t \Vert F(\tau) \Vert_{H^s} d\tau \Bigr).
\]
\end{lemma}
\begin{proof}
Let $u(t)=\tilde{S}(t,\tau)f$ be the solution at time $t$ to
\[ i\partial_t u  + \Delta u + i v \cdot \nabla u = 0,
  \quad (t, x)\in (0,T)\times \R^2 \]
with initial data $u(\tau)=f$. It suffices to show that $\tilde{S}(t,\tau)$ is bounded
on $H^s$, $0\le s \le 1$. Interpolation reduces the claim to $s=0$ and $s=1$.
If $s=0$ then
\[
\begin{split}
\frac{d}{dt} \Vert u(t) \Vert_{L^2}^2 =&  \int v \cdot \nabla |u|^2 dx
\\ \le & \Vert \nabla v (t) \Vert_{L^\infty} \Vert u(t) \Vert_{L^2}^2
\end{split}
\]
and the assertion follows by application of Gronwall's Lemma.
If $s=1$ then
\[
\begin{split}
\frac{d}{dt} \Vert  \nabla u (t) \Vert_{L^2}^2 =&
  \int v \cdot (\nabla u \Delta \overline{u} + \nabla \overline{u} \Delta u )   dx
\\ \le & 4 \Vert \nabla v (t) \Vert_{L^\infty} \Vert \nabla u(t)
\Vert_{L^2}^2
\end{split}
\]
and the assertion follows as above.
\end{proof}

We also consider the dual problem:

\[ i\partial_t u + \Delta u + i\nabla \cdot(v u) = F. \]

\begin{lemma} \label{dual} Suppose that $v$ is real valued and
\[ \int_0^T \Vert \nabla v(\tau) \Vert_{L^{\infty}} d\tau < \infty.  \]
Then, if $-1 \le s \le   0 $
\[
 \Vert u(t) \Vert_{H^s(\R)}
\le e^{ 4\int_0^t \Vert \nabla v \Vert_{L^{\infty}} d\tau }
\Bigl(\Vert u(0)
\Vert_{H^s(\R)}  + \int_0^t \Vert F(\tau) \Vert_{H^s} d\tau\Bigr).
\]
\end{lemma}
\begin{proof}
It suffices to study the case $F\equiv0$, since the general case follows by
variation of constants as above. Let $S(\tau,t)g$ be the solution to
$$
  i\partial_t u + \Delta u + i\nabla \cdot(v u) = 0, 
  \quad (\tau, x)\in (0,T)\times \R^2
$$
evaluated at time $\tau$ with initial data $u(t)=g$. We have to show
that
\begin{equation}
\label{duale}
 \Vert S(\tau,t) g \Vert_{H^{s}} \le e^{4 \int_{\tau}^t \Vert
\nabla v \Vert_{L^\infty}d\tau } \Vert g \Vert_{H^s}.
\end{equation}
Let $\tilde S(\tau,t)f$ be the solution to
\[ i\d_t u + \Delta u +  i v \cdot \nabla  u = 0,   
  \quad (t, x)\in (0,T)\times \R^2, \]
at time $t$ with initial data $u(\tau)= f$. Then $ \tilde{S}(t,\tau)$ is the adjoint
operator of $S(\tau,t)$ since
\[
\begin{split}
 \frac{d}{dt'} \langle S(t',t) f, \tilde{S}(t', \tau) g \rangle
= &  \langle i \Delta S(t',t) f, \tilde{S}(t',\tau) g\rangle
\\ &  + \langle S(t',t) f, i \Delta \tilde{S}(t',\tau) g\rangle
\\ & - \langle  \nabla\cdot \bigl(v S(t',t) f\bigr), \tilde{S}(t',\tau) g \rangle
\\ & - \langle  S(t',t) f, v\cdot \nabla\tilde{S}(t',\tau) g \rangle
\\  = & 0,
\end{split}
\]
where $\langle \cdot, \cdot \rangle$ denotes the inner product in $L^2$.
Now Lemma \ref{l2e} can be applied to $\tilde S$ and by duality we obtain
\eqref{duale}.
\end{proof}

\section{The difference of two solutions}
In this section we give the proof of the estimate \eqref{diff}.
This estimate combined with Theorem \ref{thm4.1} completes the proof
of Theorem \ref{main}.

Let
$u_0, v_0 \in H^s(\R^2)$ with $s>3/4$. According to Theorem \ref{thm4.1}
there exist solutions in
\[ L^{\infty}(0,T; H^s(\R^2)) \cap L^4(0,T; B^{1/2}_{q,2}) \]
for some $0<T\leq 1$, $q>4$.
Any solution satisfying the assumptions of Theorem
\ref{main} lies in that space. Let $u$ and $v$ be two such solutions
and let $w$ be their difference. It satisfies
\begin{equation}\label{diffeq}
\begin{split}
i\d_t w_1 + \Delta w_1 \hspace{-1.5cm}& \hspace{1.5cm}
  + 2i \nabla \cdot (A[u,u]  w_1) = -2i A[u+v,w]\cdot\nabla v_1 \\
&+ (A_0[u,u]+ |A[u,u]|^2) w_1
  + (A_0[u+v,w]+|A[u,u]|^2-|A[v,v]|^2) v_1  \\
& + 4 i \bigl\{ \Im(u_2 \overline{u}_1) w_2
+ \Im\bigl((u_2+v_2)\overline{w}_1+
w_2 (\overline{u}_1+\overline{v}_1)\bigr)v_2 \bigr\}, \\
i\d_t w_2 + \Delta w_2  \hspace{-1.5cm}& \hspace{1.5cm}
  + 2i \nabla \cdot (A[u,u]  w_2)  = -2i A[u+v,w]\cdot\nabla v_2 \\
&+ (A_0[u,u]+ |A[u,u]|^2) w_2
  + (A_0[u+v,w]+|A[u,u]|^2-|A[v,v]|^2)v_2  \\
&+  4 i \bigl\{ \Im(u_1 \overline{u}_2) w_1
  +  \Im\bigl((u_1+v_1)\overline{w}_1
  + w_1(\overline{u}_2+\overline{v}_2)\bigr)v_1 \bigr\}.
\end{split}
\end{equation}
Here, we notice that $A[u,v]$ and $A_0[u,v]$ are bilinear and symmetric
in $u$, $v$.

We want to estimate $\Vert w(t) \Vert_{H^{-1/2}}$.
Using  Lemma \ref{dual} for the equations on $w$ above, we shall show
for $t>s$ with $c$ depending on $\Vert u \Vert_{L^4 B^{1/2}_{q,2}}$,
$\Vert v \Vert_{L^4 B^{1/2}_{q,2}}$, $\Vert u \Vert_{L^\infty H^{1/2}}$, and
$\Vert v \Vert_{L^\infty H^{1/2}} $,
\begin{equation} \label{est}
\begin{split}
\Vert w(t) \Vert_{H^{-1/2}(\R^2)} \le c \Vert w(s) \Vert_{H^{-1/2}(\R^2)}
\\
&\hspace{-5cm}  + c \int_s^t (1+\Vert u(\tau)  \Vert_{B^{1/2}_{q,2}}^2 +
\Vert v(\tau) \Vert_{B^{1/2}_{q,2}}^2 )^2 \Vert w (\tau) \Vert_{H^{-1/2}} d\tau
\end{split}
\end{equation}
and hence, by Gronwall's inequality
\begin{equation}
\Vert w(t) \Vert_{H^{-1/2}}
\le   c \Vert w(0) \Vert_{H^{-1/2}}
\exp\bigl\{
(1+ \Vert u \Vert_{L^4([0,T],B^{1/2}_{q,2})}^2 +
\Vert v \Vert_{L^4([0,T],B^{1/2}_{q,2})}^2 )^2 \bigr\},
\end{equation}
which implies uniqueness of the solutions.

To establish \eqref{est} we use Lemma \ref{dual}.
We first observe that the assumption of Lemma \ref{dual}
\[
\int_0^T \Vert \nabla A[u,u](\tau) \Vert_{\infty} d\tau
\lesssim  \Vert u \Vert_{L^{4}_{T}B^{1/2}_{q,2}}^2
+   \Vert u \Vert_{L^{\infty}_{T}L^2}^2 <\infty
\]
is verified by \eqref{A}.

Then it suffices to estimate  each term on the right hand sides
of the equations \eqref{diffeq}.

The following terms are easily estimated.
\begin{equation}
\begin{split}
\Vert \Im (f \overline g ) h \Vert_{H^{-1/2}}
\lesssim & \Vert  \Im (f \overline g) \Vert_{B^{1/2}_{q,2}}
\Vert h \Vert_{H^{-1/2}} \\
\lesssim &
\Vert f \Vert_{B^{1/2}_{q,2}}
\Vert g \Vert_{B^{1/2}_{q,2}}
\Vert h \Vert_{H^{-1/2}},
\end{split}
\end{equation}
where the first inequality is a consequence of \eqref{cal3}, and the second
is a consequence of
\eqref{cal2}.
Similarly
\begin{equation}
\begin{split}
\Vert  \Im (f \overline g) h \Vert_{H^{-1/2}}
\lesssim & \Vert  \Im (f \overline g) \Vert_{H^{-1/2} }
\Vert h \Vert_{B^{1/2}_{q,2}} \\
\lesssim &
\Vert f \Vert_{B^{1/2}_{q,2}}
\Vert h \Vert_{B^{1/2}_{q,2}}
\Vert g \Vert_{H^{-1/2}},
\end{split}
\end{equation}
and using the boundedness of the Riesz transforms
\begin{equation}
\begin{split}
\Vert R_j R_k \Re (f \overline g) h \Vert_{H^{-1/2}}
\lesssim & \Vert R_j R_k \Re (f \overline g) \Vert_{B^{1/2}_{q,2}}
\Vert h \Vert_{H^{-1/2}} \\
\lesssim &
\Vert f \Vert_{B^{1/2}_{q,2}}
\Vert g \Vert_{B^{1/2}_{q,2}}
\Vert h \Vert_{H^{-1/2}}
\end{split}
\end{equation}
and
\begin{equation}
\begin{split}
\Vert R_j R_k \Re (f \overline g) h \Vert_{H^{-1/2}}
\lesssim & \Vert R_j R_k \Im (f \overline g) \Vert_{H^{-1/2} }
\Vert h \Vert_{B^{1/2}_{q,2}} \\
\lesssim &
\Vert f \Vert_{B^{1/2}_{q,2}}
\Vert h \Vert_{B^{1/2}_{q,2}}
\Vert g \Vert_{H^{-1/2}}.
\end{split}
\end{equation}

It remains to control the terms containing $A$.
In particular the estimate of the term containing  derivatives of $v$ is crucial.
This  term is controlled by Lemma \ref{Abi}:
\begin{equation}
 \| A [u,w] \nabla v  \|_{H^{-1/2}}
   \lesssim  \bigl( \| u \|_{H^{1/2}} \| v \|_{H^{1/2}}
    + \| u \|_{B^{1/2}_{q,2}} \| v \|_{B^{1/2}_{q,2}} \bigr)
   \| w \|_{H^{-1/2}}.
\end{equation}
Finally, we observe that
\[ |A[u,u]|^2 - |A[v,v]|^2 = (A[u,u]+ A[v,v]) A[u+v,w] \]
and the following two estimates complete the proof.
The first one follows from \eqref{cal3} and \eqref{Abesov},
\begin{equation}
\begin{split}
& \Vert A[u,u] A[f,g] h \Vert_{H^{-1/2}}\\
& \lesssim  \Vert A[u,u] \Vert_{B^{1/2}_{q,2}}
  \Vert A[f,g]  h \Vert_{H^{-1/2}} \\
& \lesssim  \Vert A[u,u]\Vert_{B^{1/2}_{q,2}}
 \Vert  A[f,g] \Vert_{B^{1/2}_{q,2}}
\Vert h \Vert_{H^{-1/2}}\\
& \lesssim
\bigl( \Vert u \Vert_{B^{1/2}_{q,2}}^2 + \Vert u \Vert_{L^2}^2 \bigr)
\bigl( \Vert f \Vert_{B^{1/2}_{q,2}} \Vert g \Vert_{B^{1/2}_{q,2}}
  + \Vert f \Vert_{L^2} \Vert g \Vert_{L^2} \bigr)
  \Vert h \Vert_{H^{-1/2}}.
\end{split}
\end{equation}
The second one follows from \eqref{cal3}, \eqref{Abesov}, and  Lemma \ref{Abi},
\begin{equation}
\begin{split}
& \Vert A[u,u] A[f,g] h \Vert_{H^{-1/2}}\\
& \lesssim  \Vert A[u,u] \Vert_{B^{1/2}_{q,2}}
  \Vert A[f,g] h \Vert_{H^{-1/2}} \\
& \lesssim  \Vert A[u,u] \Vert_{B^{1/2}_{q,2}} \Vert  A[f,g] h \Vert_{H^{1/2}} \\
& \lesssim \bigl( \Vert u \Vert_{B^{1/2}_{q,2}}^2 + \Vert u \Vert_{L^2}^2 \bigr)
  \bigl( \Vert f \Vert_{B^{1/2}_{q,2}} \Vert h \Vert_{B^{1/2}_{q,2}}
  + \Vert f \Vert_{H^{1/2}} \Vert h \Vert_{H^{1/2}} \bigr)
    \Vert g \Vert_{H^{-1/2}}.
\end{split}
\end{equation}
This completes the proof of \eqref{est}.

\appendix
\section{Appendix}

In this appendix, we briefly describe the derivation of the modified
Schr\"odinger map (\ref{mse}), which is due to Nahmod, Stefanov and
Uhlenbeck \cite{MR1929444}.
Recall that the Schr\"odinger map $z$ from $\R\times \R^2$ to
$(\mathbb{C}, g\, dz\, d\overline{z})\simeq \S^2$ is given by
\begin{equation}\label{sm_a}
  \frac{\d z}{\d t} =  i \sum_{j=1}^{2}
  \left( \frac{\d}{\d x_j} - 2\frac{\overline{z}\, \d_{x_j}
      z}{ 1+|z|^2} \right) \frac{\d z}{\d x_j}.
\end{equation}
It is not hard to check that
\begin{equation}\label{cons1}
  \left( \frac{\d}{\d x_j} - 2\frac{\overline{z}\, \d_{x_j}z}
{ 1+|z|^2} \right) \frac{\d z}{\d x_k} =
\left( \frac{\d}{\d x_k} - 2\frac{\overline{z}\, \d_{x_k}z}
{ 1+|z|^2} \right) \frac{\d z}{\d x_j}
\end{equation}
and
\begin{equation}\label{cons2}
 \left[  \frac{\d}{\d x_j} - 2\frac{\overline{z}\, \d_{x_j}z}{ 1+|z|^2},
 \frac{\d}{\d x_k} - 2\frac{\overline{z}\, \d_{x_k}
      z}{ 1+|z|^2}\right] = -4 i \Im (\overline{b}_j b_k)
\end{equation}
with
\[ b_j = \frac{\d_{x_j} z}{1+|z|^2} \]
hold for $j$, $k=0,1,2$.
Then we set
\[ u_j = e^{i\psi} \frac{\d_{x_j}  z}{1+|z|^2}, \]
and
\[  D_j  =  (1+|z|^2)^{-1}  e^{-i\psi} \circ \nabla_j
\circ (1+|z|^2)e^{i\psi}
   = \d_{x_j} + i A_j
\]
for $j=0,1,2$, where $x_0=t$, where the real-valued function $\psi$ is
determined later.
Note that
\[ A_j = - \d_{x_j} \psi -i \frac{z\,\d_{x_j} \overline{z}  -
  \overline{z}  \d_{x_j}  z}{1+|z|^2}= - \d_{x_j} \psi+ 2 \Im (\overline{b}_j z )
    \]
is real-valued.
By using the notation above, the equations (\ref{sm_a}), (\ref{cons1}),
and (\ref{cons2}) are rewritten as
\begin{gather}
  u_0 = i \sum_{j=1}^{2} D_j u_j \label{ms+} \\
  D_j u_k = D_k u_j, \label{cons1+} \\
  [D_j,D_k] = i \left(\frac{\d A_k}{\d x_j}
  - \frac{\d A_j}{\d x_k} \right)= -4i \Im  (u_j \overline u_k ) \label{cons2+}
\end{gather}
for $j$, $k=0,1,2$. Here, we notice that the equations
(\ref{ms+}), (\ref{cons1+}), and (\ref{cons2+}) are invariant for
arbitrary choice of $\psi$. Then, the system of the nonlinear
Schr\"odinger equations on $u_j$'s is derived as follows. For
$j=1,2$ we first notice that
$$
  D_j u_0 = D_0 u_j = \d_t u_j + i A_0 u_j
$$
holds by (\ref{cons1+}).
On the other hand, by using (\ref{ms+}), (\ref{cons1+}), and (\ref{cons2+})
we have
\[
\begin{split}
 D_j u_0
     & =  i \sum_{k=1}^2 D_j D_k u_k \\
     & = i \sum_{k=1}^{2} \bigl( D_k D_j u_k + [D_j, D_k] u_k\bigr) \\
     & =  i \sum_{k=1}^2 D_k^2 u_j + 4 \sum_{k=1}^2 \Im (u_j \overline{u}_k) u_k.
\end{split}
\]
Thus, for $j=1,2$ we obtain in a somewhat brief notation
\begin{equation}\label{msm}
  \d_t u_j = i (\nabla+ i A)^2 u_j -i A_0 u_j
    + 4 \sum_{k=1}^2 \Im (u_j \overline{u}_k) u_k,
\end{equation}
where we denote $A=(A_1, A_2)$.

Now we determine the gauge $\psi$.
For each $t$ we define $\psi$  -- up to constants -- by
\[ \Delta \psi = -i \sum_{j=1}^2 \d_{x_j} \frac{ z\,\d_{x_j} \overline{z}  -
  \overline{z}  \d_{x_j}  z}{1+|z|^2}  =  2  \sum_{j=1}^2 \d_{x_j}
  \Im( \overline{b_j} z)
\]
so that
\begin{equation}\label{coulomb}
 \nabla \cdot A=0, \qquad A_j \to 0\ \text{ as } x \to \infty.
\end{equation}
This condition and (\ref{cons2+}) enable us to determine $A$ and $A_0$ in
terms of $u_1$, $u_2$.
In fact,
\begin{equation}
  -\Delta A_1
    = \d_{x_2} \bigl(  \d_{x_1} A_2 - \d_{x_2} A_1 \bigr)
    = - 4\, \d_{x_2} \Im ( \bar{u}_1 u_2)
\end{equation}
holds by using  the first equality of (\ref{coulomb}) and  (\ref{cons2+}).
Thus,
\begin{equation}\label{a_1}
  A_1 =  4\, G_1 \ast \Im ( \bar{u}_1 u_2),
\end{equation}
where $2\pi\, G_1(x)= x_2/|x|^{2}$.
Similarly, we obtain
\begin{equation}\label{a_2}
  A_2 = 4\, G_2 \ast \Im ( \bar{u}_1 u_2)
\end{equation}
with $2\pi\, G_2(x)= - x_1/|x|^{2}$.
We use (\ref{cons2+}) and (\ref{coulomb}) again with $j=0$ to determine $A_0$,
\begin{equation}
\begin{split}
  -\Delta A_0 & = - \sum_{k=1}^2 \d_{x_k}^2 A_0\\
  & = \sum_{k=1}^2 \d_{x_k} \bigl( -\d_t A_k
    + 4 \Im( \overline{u}_k u_0) \bigr) \\
  & = 4 \sum_{k=1}^{2} \d_{x_k} \Im (\overline{u}_{k} u_0).
\end{split}
\end{equation}
Then we apply (\ref{ms+}) to obtain
\begin{align*}
  \Im(\overline{u}_{k} u_0)
  & = \Im\Bigl( \overline{u}_{k}\, i \sum_{j=1}^2 D_j u_j \Bigr)\\
  & = \sum_{j=1}^2 \Re\bigl( \overline{u}_{k} D_j u_j \bigr)\\
  & = \sum_{j=1}^{2} \bigl\{ \d_{x_j} \Re( \overline{u}_{k} u_j )
    - \Re\bigl( (\overline{D_j u_k}) u_j \bigr) \bigr\},
\end{align*}
where we used the relation $D_j f\cdot \overline{g} = \d_{x_j}(f\cdot \overline{g})
 -f\cdot \overline{D_j g}$.
Since this relation also implies
$$
  \Re\bigl( (\overline{D_j u_k}) u_j \bigr)
  = \Re\bigl( (\overline{D_k u_j}) u_j \bigr)
  = \d_{x_k} |u_j|^2/2,
$$
we obtain
$$
  -\Delta A_0 = 4 \sum_{j,k=1}^2 \d_{x_j}\d_{x_k} \Re( u_j \overline{u}_{k})
    - 2 \Delta |u|^2,
$$
where we denote $u=(u_1, u_2)$ and $|u|^2=|u_1|^2+|u_2|^2$.
Therefore,
\begin{equation}\label{a_0}
  A_0 = 4 \sum_{j,k=1}^2 R_j R_k \Re( u_j \overline{u}_{k} ) + 2 |u|^2,
\end{equation}
where $R_j$ denotes the Riesz transforms.
Therefore, we derive the system (\ref{msm}) with (\ref{a_1}), (\ref{a_2}), and
(\ref{a_0}), which is the modified Schr\"odinger map (\ref{mse}).


\begin{thebibliography}{10}


\bibitem{MR1941660} Manoussos G. Grillakis, Vagelis Stefanopoulos.
\newblock  Lagrangian formulation, energy estimates, and the Sch\"odinger map problem.
\newblock  {\em Comm. Partial Differential Equations} 27(9-10): 1845--1877, 2002.



\bibitem{MR1731461} Hiroyuki Chihara. 
\newblock  Gain of regularity for semilinear Schr\"odinger equations.
\newblock {\em Math. Ann.} 315(4): 529--567, 1999.







\bibitem{MR1737504} Nai-Heng Chang, Jalal Shatah, and Karen Uhlenbeck.
\newblock Schr\"odinger maps.
\newblock {\em Comm. Pure Appl. Math.}, 53(5):590--602, 2000.







\bibitem{MR1957040} Weiyue Ding.
\newblock On the {S}chr\"odinger flows.
\newblock In {\em Proceedings of the International Congress of Mathematicians,
  Vol. II (Beijing, 2002)}, pages 283--291, Beijing, 2002. Higher Ed. Press.







\bibitem{MR1877231} Weiyue Ding and Youde Wang.
\newblock Local {S}chr\"odinger flow into {K}\"ahler manifolds.
\newblock {\em Sci. China Ser. A}, 44(11):1446--1464, 2001.







\bibitem{MR1255899} Nakao Hayashi and Tohru Ozawa.
\newblock  Remarks on nonlinear {S}chr\"odinger equations in one space dimension.
\newblock {\em Differential Integral Equations} 7(2): 453--461, 1994.







\bibitem{0920.58022} Fr{\'e}d{\'e}ric H{\'e}lein.
\newblock {\em {Applications harmoniques, lois de conservation et rep{\`e}res
  mobiles. Pr{\'e}face de James Eells. (Harmonic mappings, conservation laws
  and moving frames. Preface by James Eells).}}
\newblock {Nouveaux Essais. Paris: Diderot Editeur. xix, 286 p. FF 180.00},
  1996.







\bibitem{Kato} Jun Kato.
\newblock Existence and uniqueness of the solution to the modified
  Schr\"{o}dinger map.
\newblock {\em Math. Res. Lett.}, to appear, 2005.







\bibitem{NK} Carlos E.~Kenig, Andrea Nahmod.
\newblock Personal communication.
\newblock 2004.







\bibitem{MR1230709} Carlos E.~Kenig, Gustavo Ponce, and Luis Vega.
\newblock  Small solutions to nonlinear Schr\"odinger equations.
\newblock {\em Ann. Inst. H. Poincar\'e Anal. Non Lin\'eaire}, 10(3): 255--288, 1993.







\bibitem{MR2094472} Joachim Krieger.
\newblock Global regularity of wave maps from {$\mathbf{R}\sp {2+1}$} to {$H\sp 2$}. {S}mall energy.
\newblock {\em Comm. Math. Phys.}, 250(3):507--580, 2004.







\bibitem{MR1929444} Andrea Nahmod, Atanas Stefanov, and Karen Uhlenbeck.
\newblock On {S}chr\"odinger maps.
\newblock {\em Comm. Pure Appl. Math.}, 56(1):114--151, 2003.







\bibitem{MR2038118} Andrea Nahmod, Atanas Stefanov, and Karen Uhlenbeck.
\newblock Erratum: ``{O}n {S}chr\"odinger maps'' [{C}omm. {P}ure {A}ppl.\
  {M}ath.\ {\bf 56} (2003), no.\ 1, 114--151; mr 1929444].\newblock {\em Comm. Pure Appl. Math.}, 57(6):833--839, 2004.







\bibitem{MR1696311} Catherine Sulem and Pierre-Louis Sulem.
\newblock {\em The nonlinear {S}chr\"odinger equation}, volume 139 of {\em
  Applied Mathematical Sciences}.
\newblock Springer-Verlag, New York, 1999.







\bibitem{MR866199} Pierre-Louis Sulem, Catherine Sulem, and Claude Bardos.
\newblock On the continuous limit for a system of classical spins.
\newblock {\em Comm. Math. Phys.}, 107(3):431--454, 1986.


\bibitem{MR1869874}
Terence Tao.
\newblock Global regularity of wave maps. {II}. {S}mall energy in two
  dimensions.
\newblock {\em Comm. Math. Phys.}, 224(2):443--544, 2001.

\bibitem{MR2043751}
Daniel Tataru.
\newblock The wave maps equation.
\newblock {\em Bull. Amer. Math. Soc. (N.S.)}, 41(2):185--204 (electronic),
  2004.



\end{thebibliography}
 \end{document}